\documentclass[12pt]{article}
\usepackage{amssymb,amsmath}
\pdfoutput=1

\usepackage[pdftex]{graphicx}         

\newcommand{\pr}{\mathbb{P}}                     
\newcommand{\exn}{\mathbb{E}\,}                  
\newcommand{\var}{\mathbb{V}\,}
\newcommand{\deq}{\stackrel{d}{=}}               

\newtheorem{theorem}{Theorem}

\newcommand{\proof}[1]{{{\em Proof}#1.}}
\newcommand{\eproof}{\hfill$\Box$\par\medskip\smallskip}

\begin{document}
\title{On the  ruin time distribution for a Sparre Andersen process with exponential claim sizes}
\author{Konstantin A Borovkov\footnote{Research supported by the ARC Centre
of Excellence for Mathematics and Statistics of Complex Systems.} \ and David C M
Dickson}
\date{}
\maketitle

\begin{abstract}
We derive a closed-form (infinite series) representation for the distribution of the
ruin time for the Sparre Andersen model with exponentially distributed claims. This
extends a recent result of Dickson et al.~\cite{DiHuZh05} for such processes with Erlang
inter-claim times. We illustrate our result in the cases of gamma and mixed exponential
inter-claim time distributions.
\end{abstract}



\section{Introduction}

In the Sparre Anderson model, the (continuous-time) surplus process $\{U(t)\}_{t\ge 0}$
has the form
\[
U (t) = u + ct - \sum_{j\le N(t)} X_j,
\]
where $u\ge0$ is the initial surplus, $c>0$ is the premium rate, and  $\{N(t)\}_{t\ge 0}$ is a
delayed renewal process generated by a sequence of inter-claim times $\{T_j\}_{j\ge 0}$:
\[
N(t) = \inf\{ j \ge 0: \, T_0 +\cdots + T_j \ge t\},
\]
and $\{X_j\}_{j\ge 1}$ is the sequence of claim sizes (so that a claim of size $X_1$ is
made at time $T_0$, etc). We assume that the random variables from the above sequences
are jointly independent, with $\{T_j\}_{j\ge 1}$ and $\{X_j\}_{j\ge 1}$ being
i.i.d.\ sequences. The goal of the present note is to derive an explicit formula for the
distribution of the ruin time
\[
\tau = \inf\{ t>0: \, U(t) <0\}
\]
in the special case when the $X_j$'s follow the exponential distribution.

When claims occur according to a Poisson process and the claim size distribution is exponential, a solution for the distribution of the ruin time $\tau$
has been known for many years. See, for example, \cite{Asm00}, \cite{DreWil03} and \cite{Seal78} for different solutions to this problem.
In the recent paper \cite{DiHuZh05}, the authors
used analytical techniques to obtain an explicit formula for the density of $\tau$ in
the case when the $T_j$'s have an Erlang distribution. In the present note, we present
an alternative probabilistic method, which enables one to derive such an explicit
formula in the more general case when the $T_j$'s follow an arbitrary distribution.

\section{The main result}

We assume that the claim sizes $X_j$ follow the exponential distribution with
parameter~$\lambda>0$:
\begin{equation}
\pr (X_j > x)=e^{-\lambda x },\qquad x\ge 0,
\label{ExpX}
\end{equation}
while the positive random variables $T_0$ and $T_1 (\deq T_j$, $j>1)$   have densities
$f_0(t)$ and $f(t)$, respectively. By $g*h$ we denote the convolution of the functions
$g,h$ defined on $(0,\infty)$:
\[
(g*h) (t) = \int_0^t g(t-v) h(v) dv,
\]
and by $g^{*n}= g^{*(n-1)} *g,$ $n\ge 2$, the $n$-fold convolution of $g$ with itself.

\begin{theorem}
Under the above assumptions, the ruin time $\tau$ has a (defective) density $p_\tau (t)$ given
by
\begin{multline}
p_\tau (t) = e^{-\lambda (u+ct)} \biggl\{ f_0(t)
 \\
 +
 \sum_{n=1}^\infty \frac{\lambda^n (u+ct)^{n-1}}{n!} \bigl[ u (f^{*n} * f_0 )(t)
   + c(f^{*n} * f_1 ) (t)\bigr]\biggr\},
\label{Main}
\end{multline}
where $f_1(t) = tf_0(t).$
\end{theorem}

\medskip

\proof{} The idea of the proof is similar to the one used in~\cite{Bo85}: first we will
translate our problem into the problem of the crossing of a linear boundary by the pure jump process
$U^0 (t) =U(t)-ct$ and then swap the roles of the time and space coordinates. Then we
notice that the generalised inverse of the function $U^0(t)$ is nothing else but the
trajectory of a compound Poisson process. Eventually, the original problem proves to be
equivalent to finding the distribution of the hitting time of a level by a skip-free
L\'evy process, of which the solution is well-known and is given by Kendall's identity (see
e.g. \S\,12, Theorem~1 in~\cite{Bo72}, or \cite{BoBu01}).

\medskip

(i)~We will assume in parts (i)-(ii) of the proof that $T_0\equiv v=\text{const}$ (which
is equivalent to conditioning on $T_0$, but is more convenient from a notation
viewpoint).

As we have just said, it is easily seen that, for the pure jump process
\[
U^0 (t) =U(t)-ct \equiv u - \sum_{j\le N(t)} X_j,
\]
one has
\[
\tau = \inf\{ t>0: \, U^0 (t) - (-ct) <0\}.
\]
Next we `translate' the origin to the point $(v,u)$ and swap the roles of coordinates by
introducing the new `time' $s=u-x$ and `space' $y=t-v$ (where $t$ and $x$ respectively represent the original time and space). In the new system of
coordinates, the trajectory of our process $\{U^0 (t)\}$ is again a pure step function,
which starts at zero at `time' $s=0$ and has jumps of sizes $T_1, T_2,T_3,\dots,$ at
`times' $X_1, X_1+ X_2, X_1+ X_2 +X_3,\dots.$ Due to our assumption~\eqref{ExpX}, this
will be a trajectory of the compound Poisson process
\[
Z^0 (s) = \sum_{k\le M(s)}  T_k,
\]
where
\[
M(s) = \inf\{ k\ge 1:\, X_1 +\cdots +X_k > s\} -1
\]
is a Poisson process with rate $\lambda$. The distribution of
the r.v.\ $Z^0 (s)$ with $s>0$ has an atom $e^{-\lambda s}$ at zero and a density on $(0,\infty)$ given by
\begin{equation}
p_{Z^0 (s)} (y)= e^{-\lambda s}\sum_{n=1}^\infty \frac{(\lambda s)^n}{n!} f^{*n}(y),
\qquad y>0.
\label{DenZ0}
\end{equation}

To a crossing of the (lower) linear boundary $x=-ct$ by the process $\{U^0 (t)\}$ at time
$\tau$ (this necessarily is a jump epoch) there corresponds a (continuous) crossing of the
(again lower linear) boundary
\[
y = s/c - (v+u/c), \qquad s>0,
\]
by the process $\{Z^0 (s)\}$ at `time' $\sigma = u+c\tau$, so that
\begin{equation}
\tau = (\sigma -u)/c.
\label{TauSigma}
\end{equation}

Finally, we notice that $\sigma$ is the crossing time of the (lower) level   $-(v+u/c)$
by the process $Z  (s)= Z^0 (s) -s/c,$ which is clearly a skip-free in the negative
direction L\'evy process.

Figure 1 illustrates the translation of the original problem. The original surplus process starts at level $u=1$, and ruin occurs at the fifth claim.
Rotating the figure anti-clockwise through 90 degrees we see the corresponding path of the pure jump process $\{Z^0(s)\}$.

\begin{figure}[h]
\begin{center}
%
%
%
%
\mbox{\pdfimage width 0.95\linewidth {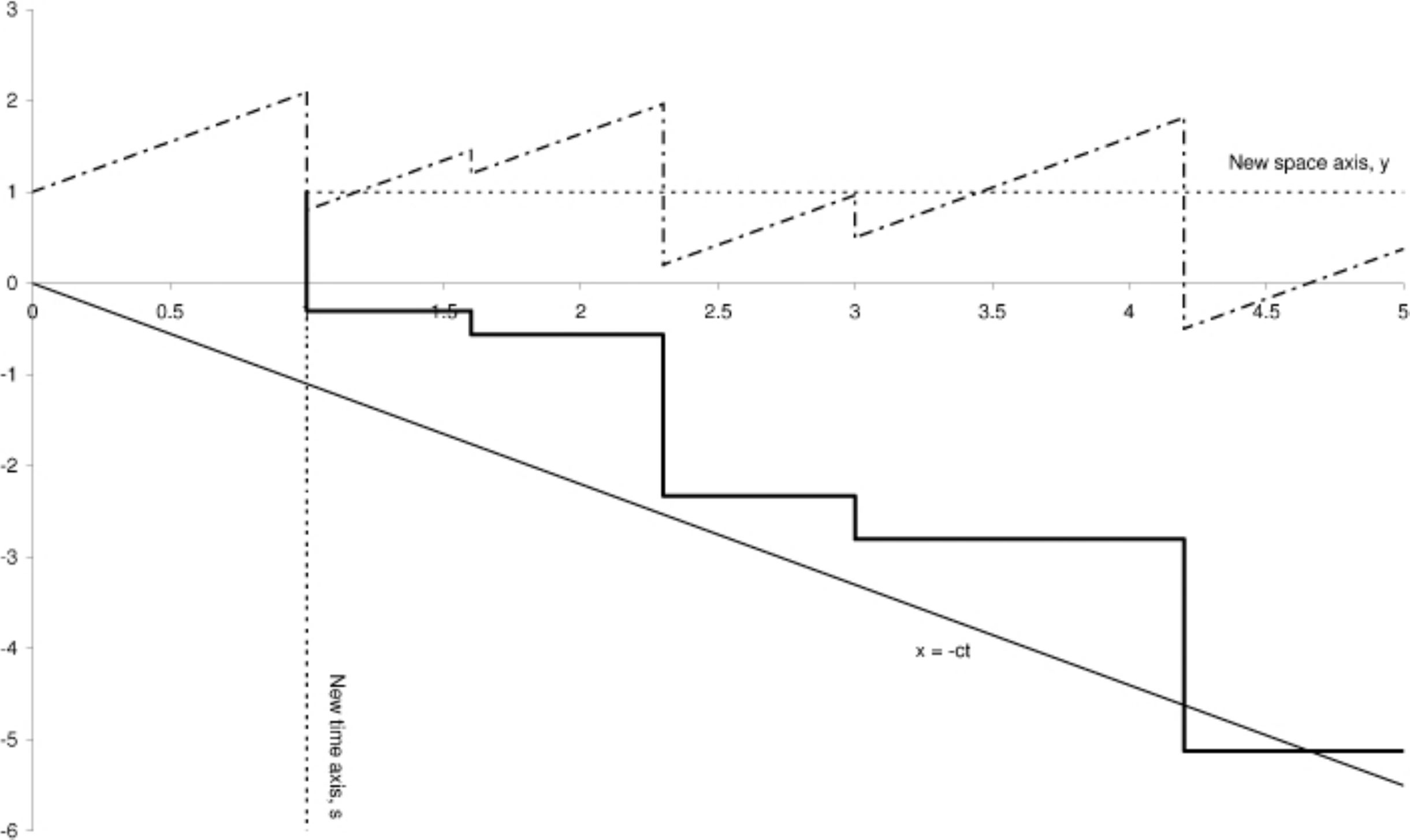}}
\caption{Original and translated processes}
\end{center}
\label{Fig1}
\end{figure}

\medskip

(ii)~Therefore, provided that $Z(s)$ has a density $p_{Z(s)} (y)$ at the point
$y=-(v+u/c),$ the crossing `time' $\sigma$ also has a density $p_\sigma (s)$ at the
point $s$, which is given by Kendall's identity (see e.g. \S\,12, Theorem~1
in~\cite{Bo72}, or \cite{BoBu01}):
\[
p_\sigma (s) = \frac{v+u/c}{s} p_{Z(s)} (-(v+u/c)).
\]
This together with \eqref{DenZ0} implies that, for $s> u+cv,$ the r.v. $\sigma$ has the
density
\[
p_\sigma (s) = \frac{v+u/c}{s}\, e^{-\lambda s}
 \sum_{n=1}^\infty \frac{(\lambda s)^n}{n!} f^{*n} \bigl((s-u)/c -v\bigr).
\]
Therefore, it follows now from~\eqref{TauSigma} that, given $T_0 =v$, for $t>v$ the stopping time
$\tau$ has a conditional density given by
\begin{equation}
p_\tau (t|v) = c p_\sigma (u+ct)
 = \frac{u + cv}{u+ct}\, e^{-\lambda (u+ct)}
 \sum_{n=1}^\infty \frac{(\lambda (u+ct))^n}{n!} f^{*n}  (t -v ).
\label{DenTauV}
\end{equation}

(iii) To obtain the density of $\tau$ in the general case, we observe that $\tau
\ge T_0$ always and so, using \eqref{ExpX},
\begin{multline*}
\pr (\tau \le t) = \pr (T_0 = \tau \le t) + \pr (T_0 < \tau \le t)
 \\
 = \int_0^t \pr (u+cv -X_1 <0) f_0 (v) dv
 + \int_0^t \pr (v <\tau \le t|\, T_0 =v) f_0 (v) dv
 \\
 = \int_0^t e^{-\lambda (u+cv)} f_0 (v) dv
 + \int_0^t \biggl[ \int_v^t p_\tau (r|v) dr\biggr] f_0 (v) dv.
\end{multline*}
Differentiating both sides and substituting the representation for $ p_\tau (t|v)$ from
\eqref{DenTauV}  yields the density of $\tau$:
\begin{multline*}
p_\tau (t)  =  e^{-\lambda (u+ct)} f_0 (t)
 + \int_0^t   p_\tau (t|v)   f_0 (v) dv
 \\
 =  e^{-\lambda (u+ct)} \biggl[ f_0 (t) +  \frac{1}{u+ct}
 \sum_{n=1}^\infty \frac{(\lambda (u+ct))^n}{n!} \int_0^t (u + cv) f^{*n}  (t -v ) f_0 (v) dv
 \biggr]
\end{multline*}
(the change of the order of integration/summation is justified as the integrand is a
non-negative function). As the last expression is equivalent to the RHS
of~\eqref{Main}, the theorem is proved. \eproof

\section{Examples}
\subsection {Gamma inter-claim times}
Let us first consider the situation where claims occur according to an ordinary renewal process, so that
each $T_j$, $j=0,1,2,...$ has density function
\[
f(t)=f_{0}(t)=\frac{\beta ^{n}t^{n-1}e^{-\beta t}}{\Gamma (n)},
\]%
where $n>0$ and $\beta >0$. It is well known that%
\[
f^{\ast (m+1)}(t)=\frac{\beta ^{n(m+1)}t^{n(m+1)-1}e^{-\beta t}}{\Gamma
(n(m+1))}
\]%
and it is straightforward to show that%
\[
f^{\ast m}\ast f_{1}(t)=\frac{n}{\beta }\frac{\beta
^{n(m+1)+1}t^{n(m+1)}e^{-\beta t}}{\Gamma (n(m+1)+1)}.
\]%
Then formula \eqref{Main} gives%
\begin{eqnarray*}
p_{\tau }(t)
&=&(\beta t)^{n-1}\frac{u\beta e^{-\lambda (u+ct)-\beta t}}{u+ct}%
\sum_{m=0}^{\infty }\frac{\lambda ^{m}(u+ct)^{m}}{m!}\frac{(\beta t) ^{nm}}{%
\Gamma (n(m+1))} \\
&&+ (\beta t)^{n}\frac{c n e^{-\lambda (u+ct)-\beta t}}{u+ct%
}\sum_{m=0}^{\infty }\frac{\lambda ^{m}(u+ct)^{m}}{m!}\frac{(\beta t)^{nm}%
}{\Gamma (n(m+1)+1)}.
\end{eqnarray*}
In the special case when $n$ is a positive integer, we can compute this as
\begin{multline}
p_{\tau }(t)= \frac{\beta e^{-\lambda (u+ct)-\beta t}}{u+ct}\frac{(\beta t)^{n-1}}{%
\Gamma (n)} \left( u\ _{0}F_{n}\left( 1,1+\frac{1}{n},...,1+\frac{n-1}{n};\frac{%
\lambda (u+ct)(\beta t)^{n}}{n^{n}}\right) \right.
\\
+\left.%
ct \ _{0}F_{n}\left( 1+\frac{1}{n},1+\frac{2}{n},...,1+\frac{n}{n};%
\frac{\lambda (u+ct)(\beta t)^{n}}{n^{n}}\right) \right), \label{orderl2}
\end{multline}%
where
\begin{equation*}
{}_{p}F_{q}(B_1,B_2,...,B_p,C_{1},C_{2},\ldots C_{q};Z)=\sum_{m=0}^{\infty }\frac{(B_1)_m (B_2)_m...(B_p)_m}{%
(C_{1})_{m}(C_{2})_{m}...(C_{q})_{m}}\frac{Z^{m}}{m!}
\end{equation*}%
is the generalised hypergeometric function (and $(a)_n =\Gamma(a+n)/\Gamma(a)$ is Pochhammer's symbol).
Formula \eqref{orderl2} follows from the identity
\[
\frac {\Gamma(n+1)} {\Gamma((n(m+1)+1)}=\frac{1}{n^{nm}}\prod_{k=0}^{n-1}
\frac {\Gamma(1+\frac{k+1}{n})}{\Gamma(m+1+\frac{k+1}{n})} ,
\]
which can be derived by applying the multiplication formula of Gauss as
described in~\cite{DiHuZh05}.

Formula \eqref{orderl2} is in a different form to the formula for $p_\tau (t)$ derived in~\cite{DiHuZh05}. A comparison of these
two formulae for $p_\tau (t)$ yields the identity

\begin{eqnarray*}
&&_{0}F_{n}\left( 1,1+\frac{1}{n},...,1+\frac{n-1}{n};\frac{%
\lambda (u+ct)(\beta t)^{n}}{n^{n}}\right) \\
&-& _{0}F_{n}\left( 1+\frac{1}{n},1+\frac{2}{n},...,1+\frac{n}{n};\frac{%
\lambda (u+ct)(\beta t)^{n}}{n^{n}}\right) \\
&=& \lambda (u+ct)(\beta t)^{n} \frac{n!}{(2n)!} \ _{0}F_{n}\left( 2+\frac{1}{n},2+\frac{2}{n},...,2+\frac{n}{n};\frac{%
\lambda (u+ct)(\beta t)^{n}}{n^{n}}\right).
\end{eqnarray*}
In the special case $n=1$, by writing $z=\sqrt{4\lambda\beta t (u+ct)}$ this identity reduces to the well-known result (e.g.~\cite{AbSt65})
\[
I_0 (z)-\frac{2}{z} I_1 (z)=I_2(z),
\]
where $I_v$ is the modified Bessel function of order $v$.

Next, let us consider the special case when $n=2$, and let us further assume that claims occur according to
a stationary renewal process, so that the distribution of $T_0$ is the equilibrium distribution of
$T_1$. Then we find that
\[
f_{0}(t)=\frac{\beta }{2}e^{-\beta t}\left( 1+\beta t\right) =\frac{1}{2}%
(\beta e^{-\beta t}+\beta ^{2}te^{-\beta t}),
\]%
giving
\[
f^{\ast m}\ast f_{0}(t)=\frac{1}{2}\left( \frac{\beta ^{2m+1}t^{2m}e^{-\beta
t}}{\Gamma (2m+1)}+\frac{\beta ^{2m+2}t^{2m+1}e^{-\beta t}}{\Gamma (2m+2)}%
\right) .
\]%
Further,%
\[
f_{1}(t)=tf_{0}(t)=\frac{1}{2}(\beta te^{-\beta t}+\beta ^{2}t^{2}e^{-\beta
t}),
\]%
giving%
\[
f^{\ast m}\ast f_{1}(t)=\frac{1}{2}\frac{\beta ^{2m+1}t^{2m+1}e^{-\beta t}}{%
\Gamma (2m+2)}+\frac{\beta ^{2m+2}t^{2m+2}e^{-\beta t}}{\Gamma (2m+3)}.
\]%
Then formula \eqref{Main} gives

\begin{eqnarray*}
p_{\tau }(t) &=&e^{-\lambda (u+ct)}\left( f_{0}(t)+\frac{u}{u+ct}%
\sum_{m=1}^{\infty }\frac{\lambda ^{m}(u+ct)^{m}}{m!}\frac{1}{2}\left( \frac{%
\beta ^{2m+1}t^{2m}e^{-\beta t}}{\Gamma (2m+1)}+\frac{\beta
^{2m+2}t^{2m+1}e^{-\beta t}}{\Gamma (2m+2)}\right) \right.  \\
&&\hspace{0.75in}\left. +\frac{c}{u+ct}\sum_{m=1}^{\infty }\frac{\lambda
^{m}(u+ct)^{m}}{m!}\left( \frac{1}{2}\frac{\beta ^{2m+1}t^{2m+1}e^{-\beta t}%
}{\Gamma (2m+2)}+\frac{\beta ^{2m+2}t^{2m+2}e^{-\beta t}}{\Gamma (2m+3)}%
\right) \right) ,
\end{eqnarray*}%
and we can incorporate $f_{0}(t)$ into the sums so that both start at $\ m=0$. For computational
purposes we can write this in terms of
generalised hypergeometric functions as
\begin{multline}
p_{\tau }(t) =\frac{\beta e^{-\lambda (u+ct)-\beta t}}{2(u+ct)}\left( u\
_{0}F_{2}\left( \frac{1}{2},1;\frac{\lambda (u+ct)(\beta t)^{2}}{4}\right)
\right.  \\
\left. +t(\beta u+c)\ _{0}F_{2}\left( 1,\frac{3}{2};\frac{\lambda
(u+ct)(\beta t)^{2}}{4}\right) +c\beta t^{2}\ _{0}F_{2}\left( \frac{3}{2},2;%
\frac{\lambda (u+ct)(\beta t)^{2}}{4}\right) \right).  \label{staterl2}
\end{multline}

Table \ref{table1} shows some values of finite time ruin probabilities when $\lambda=1$, $\beta=2$ and $c=1.1$. We
use the notation $\psi (u,t)$ to denote the probability of ruin by time $t$
from initial surplus $u$ when the density of $\tau $ is given by formula \eqref{orderl2}
with $n=2$, and $\psi _{e}(u,t)$ denotes the corresponding probability when
the density of $\tau $ is given by formula \eqref{staterl2}. These values have been found
by integrating the density functions numerically using Mathematica. We can
observe from this table that for each combination of $u$ and $t$, the finite
time ruin probability is greater when the distribution of $T_{0}$ is the
equilibrium distribution of $T_{1}$. This arises because both the mean and variance of $T_0$ are
smaller than the corresponding values for $T_1$.\bigskip

\begin{table}[h] \centering%
\begin{tabular}{ccccccc}
$t$ & $\psi (0,t)$ & $\psi _{e}(0,t)$ & $\psi (10,t)$ & $\psi _{e}(10,t)$ & $%
\psi (20,t)$ & $\psi _{e}(20,t)$ \\ \hline
$20$ & 0.7973 & 0.8463 & 0.0457 & 0.0509 & 0.0009 & 0.0010 \\
$40$ & 0.8332 & 0.8735 & 0.1008 & 0.1082 & 0.0060 & 0.0066 \\
$60$ & 0.8481 & 0.8848 & 0.1387 & 0.1469 & 0.0138 & 0.0148 \\
$80$ & 0.8564 & 0.8912 & 0.1651 & 0.1737 & 0.0218 & 0.0232 \\
$100$ & 0.8618 & 0.8952 & 0.1842 & 0.1930 & 0.0292 & 0.0309 \\  \hline%
\end{tabular}%
\caption{Finite time ruin probabilities.
\label{table1}}%
\end{table}%

\subsection{Mixed exponential inter-claim times}

Let us now consider the situation when the distribution of each $T_{j}$, $%
j=0,1,2,...,$ is mixed exponential with density function%
\[
f(t)=f_0(t)=p\alpha e^{-\alpha t}+q\beta e^{-\beta t},
\]%
where $0<p<1$, $q=1-p$, and $\beta >\alpha >0$. Following ideas in \cite{WilWoo07}, it is shown in \cite{Dic07} that the $m$-fold convolution of $f$ with
itself as can be written as%
\[
f^{\ast m}(t)=\sum_{j=0}^{\infty }\gamma _{m,j}\ e(m+j,\beta ;t),
\]%
where $e(m,\beta ;t)$ denotes the Erlang($m$) density with scale parameter $%
\beta $ and
\[
\gamma _{m,j}=q^{m}(1-\alpha /\beta )^{j}\sum_{r=0}^{m}\binom{m}{r}\frac{%
(r)_{j}}{j!}\left( \frac{\alpha p}{\beta q}\right) ^{r}.
\]%
We can find a similar type of expression for $f^{\ast m}\ast f_{1}(t)$ by
using Laplace transforms. For a function $w$, let%
\[
\tilde{w}(s)=\int_{0}^{\infty }e^{-st}w(t)dt.
\]%
Then%
\[
\tilde{f}(s)=\frac{p\alpha }{\alpha +s}+\frac{q\beta }{\beta +s}
\]%
and
\[
\tilde{f}_{1}(s)=\frac{p\alpha }{(\alpha +s)^{2}}+\frac{q\beta }{(\beta
+s)^{2}},
\]%
leading to%
\begin{eqnarray*}
\left[ \tilde{f}(s)\right] ^{m}\ \tilde{f}_{1}(s) &=&\frac{p}{\alpha }%
\sum_{r=0}^{m}\binom{m}{r}p^{r}q^{m-r}\left( \frac{\alpha }{\alpha +s}%
\right) ^{r+2}\left( \frac{\beta }{\beta +s}\right) ^{m-r} \\
&&+\frac{q}{\beta }\sum_{r=0}^{m}\binom{m}{r}p^{r}q^{m-r}\left( \frac{\alpha
}{\alpha +s}\right) ^{r}\left( \frac{\beta }{\beta +s}\right) ^{m-r+2}.
\end{eqnarray*}%
Hence%
\begin{eqnarray*}
f^{\ast m}\ast f_{1}(t) &=&\frac{p}{\alpha }\sum_{r=0}^{m}\binom{m}{r}%
p^{r}q^{m-r}\int_{0}^{t}\frac{\alpha ^{r+2}y^{r+1}e^{-\alpha y}}{\Gamma (r+2)%
}\frac{\beta ^{m-r}(t-y)^{m-r-1}e^{-\beta (t-y)}}{\Gamma (m-r)}dy \\
&&+\frac{q}{\beta }\sum_{r=0}^{m}\binom{m}{r}p^{r}q^{m-r}\int_{0}^{t}\frac{%
\alpha ^{r}y^{r-1}e^{-\alpha y}}{\Gamma (r)}\frac{\beta
^{m-r+2}(t-y)^{m-r+1}e^{-\beta (t-y)}}{\Gamma (m-r+2)}dy \\
&=&p\alpha \sum_{r=0}^{m}\binom{m}{r}\left( \alpha p\right) ^{r}\left( \beta
q\right) ^{m-r}\frac{e^{-\beta t}t^{m+1}}{\Gamma (m+2)}\
_{1}F_{1}(r+2,m+2,(\beta -\alpha )t) \\
&&+q\beta \sum_{r=0}^{m}\binom{m}{r}\left( \alpha p\right) ^{r}\left( \beta
q\right) ^{m-r}\frac{e^{-\beta t}t^{m+1}}{\Gamma (m+2)}\
_{1}F_{1}(r,m+2,(\beta -\alpha )t).
\end{eqnarray*}
If we now replace the $ _1 F_1$ functions by their series representations, we find after a small amount of manipulation, that
\[
f^{\ast m}\ast f_{1}(t) = \sum_{i=0}^{\infty }\eta _{i,m}e(m+i+2,\beta ;t),
\]
where%
\[
\eta _{i,m}=\frac{q^{m}}{\beta ^{2}}\left( 1-\alpha /\beta \right)
^{i}\sum_{r=0}^{m}\binom{m}{r}\frac{(r)_{i}}{i!}\left( \frac{\alpha p}{%
\beta q}\right) ^{r} \left( (r+2)(r+1)\alpha p + \beta q \right).
\]%
Thus, we have formulae for all the ingredients in formula \eqref{Main}.

Figure 2 shows three plots of the density of $\tau$ when $u=10$, $c=1.1$ and the parameters of the mixed
exponential distribution are as in Table \ref{table2}. In Figure 2, Case A is illustrated by the dotted line,
Case B by the solid line, and Case C by the bold line. We observe that the ordering of these three plots leads
to highest finite time ruin probabilities for Case A and lowest for Case C, consistent with the ordering of the
three values of $\var[T_{0}]$.

\begin{table}[h] \centering%
\begin{tabular}{cccccc}
Case & $\alpha $ & $\beta $ & $p$ & $\exn[T_{0}]$ & $\var[T_{0}]$ \\  \hline
A & 2/5 & 2 & 1/4 & 1 & 5/2 \\
B & 1/2 & 2 & 1/3 & 1 & 2 \\
C & 3/5 & 2 & 3/7 & 1 & 5/3\\  \hline
\end{tabular}
\caption{Parameters of mixed exponential distributions.
\label{table2}}%
\end{table}%

\begin{figure}[h]
\begin{center}
%
%
%
%
\mbox{\pdfimage width 1\linewidth {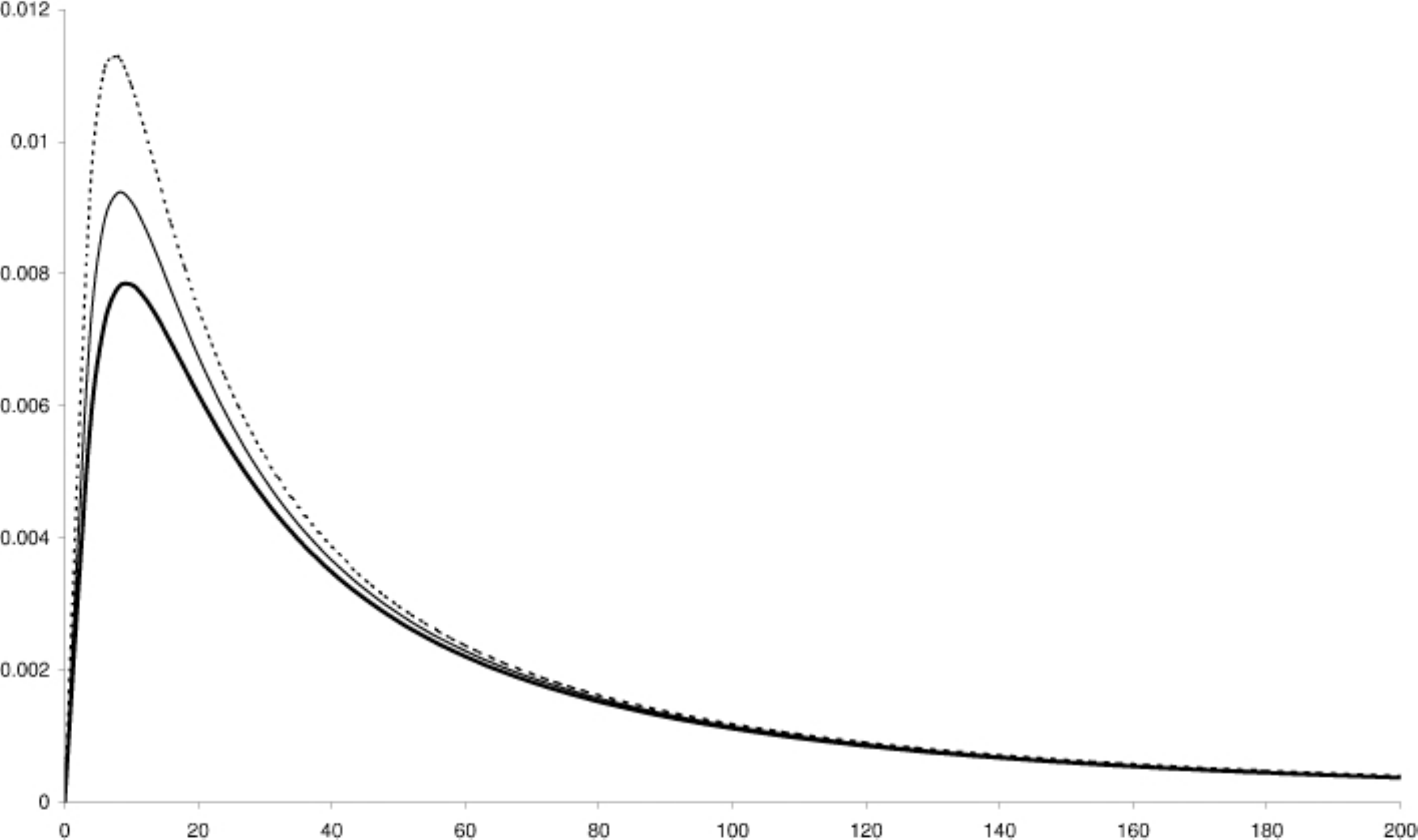}}
\caption{Densities of time to ruin for mixed exponential inter-claim times.}
\end{center}
\label{Fig2}
\end{figure}


\bigskip

\noindent Konstantin A Borovkov\newline
Department of Mathematics and Statistics\newline
The University of Melbourne\newline
Victoria 3010\newline
Australia\newline
K.Borovkov@ms.unimelb.edu.au\newline

\medskip
\noindent David C M Dickson\newline
Centre for Actuarial Studies\newline
Department of Economics\newline
The University of Melbourne\newline
Victoria 3010\newline
Australia\newline
dcmd@unimelb.edu.au

\end{document}